\documentclass[preprint, 12pt, authoyear]{elsarticle}
\usepackage{amssymb}
\usepackage[all,arc]{xy}
\usepackage{mathrsfs}
\usepackage[mathscr]{eucal}
\usepackage{amsbsy}
\usepackage{amsfonts}
\usepackage{amsmath}
\usepackage{amssymb}
\usepackage{amsthm}
\usepackage{array}
\usepackage{graphicx}

\DeclareMathOperator*{\argmin}{arg\,min}
\newcommand{\cH}{\mathcal{H}}

\newcommand{\norm}[1]{\left\Vert #1 \right\Vert}

\usepackage{bbm}
\usepackage{mathtools} 

\newtheorem{thm}{Theorem}[section]

\makeatletter
\def\ps@pprintTitle{%
 \let\@oddhead\@empty
 \let\@evenhead\@empty
 \def\@oddfoot{}%
 \let\@evenfoot\@oddfoot}
\makeatother

\usepackage{hyperref}

\begin{document}

\begin{frontmatter}
\title{A theoretical guarantee for data completion via geometric separation.}

\author[auth1]{Emily J. King}
\author[auth2]{James M. Murphy}

\address[auth1]{Center for Industrial Mathematics, University of Bremen, Bremen, Germany}

\address[auth2]{Department of Mathematics, Johns Hopkins University, Baltimore, USA}

\begin{abstract}
Scientific and commercial data is often incomplete.  Recovery of the missing information is an important pre-processing step in data analysis. Real-world data can in many cases be represented as a superposition of two or more different types of structures.  For example, images may often be decomposed into texture and cartoon-like components. When incomplete data comes from a distribution well-represented as a mixture of different structures, a sparsity-based method combining concepts from data completion and data separation can successfully recover the missing data. This short note presents a theoretical guarantee for success of the combined separation and completion approach which generalizes proofs from the distinct problems.
\end{abstract}

\end{frontmatter}

For a Hilbert space $\cH$ and coefficient index set $I$, a \emph{Parseval frame} $\Phi: \cH \rightarrow \ell^2(I)$ (or more precisely, the analysis operator of a Parseval frame) is such that $\norm{\Phi x}_2 = \norm{x}_{\cH}$ for all $x \in \cH$.  Data in a variety of applications of interest admit sparse representations with respect to a Parseval frame.  It is known, for example, that if $\Phi$ is a shearlet frame for $L^2([0,1]^2)$ and $x$ is a cartoon-like image, then there is some relatively ``small'' subset $\Lambda \subset I$ such that  $\mathbbm{1}_{\Lambda} \Phi x$  captures ``most'' of the important information of $x$, where $\mathbbm{1}_{\Lambda}$ is the function that is $1$ on $\Lambda$ and $0$ off of $\Lambda$ \cite{Guo2007}. This sparsity can be used to perform image and data processing tasks including geometric separation of data into components with different fundamental structures \cite{SPIE_IP,DK10} and image inpainting or, more broadly, data recovery \cite{SPIE_IP,JIV_IP,KKL13}.  These two tasks can be combined.  For example, by noting that natural images are usually a superposition of a cartoon-like part which is sparsely represented by a shearlet or curvelet frame and a textured part which is sparsely represented by a discrete cosine frame, one can use both frames simultaneously to perform empirically successful data completion \cite{ejkESQD05,StuerckMS}.  What follows is a generalization to the combined problem of the theoretical guarantees for data recovery \cite{SPIE_IP,JIV_IP} and for data separation \cite{SPIE_IP,DK10}.

\begin{thm}\label{thm1}Let $\mathcal{H}=\mathcal{H}_{K}\bigoplus\mathcal{H}_{M}$ be an orthogonal decomposition of a Hilbert space $\cH$. Let $\Phi_1$ and $\Phi_2$ be (the analysis operators of) Parseval frames for $\cH$ index sets $I_1$, $I_2$, respectively.  Suppose $x^{0}\in\mathcal{H}$ may be decomposed as $x^{0}=x_{1}^{0}+x_{2}^{0}$ and for $\delta\geq 0$ there exist some $\Lambda_1 \subset I_1$, $\Lambda_2\subset I_2$ such that the following $\delta$-sparsity condition is satisfied: 
\[
\label{sparsity2}\|\mathbbm{1}_{\Lambda_{1}^{c}}\Phi_{1}x_{1}^{0}\|_{1}+\|\mathbbm{1}_{\Lambda_{2}^{c}}\Phi_{2}x_{2}^{0}\|_{1}\le\delta.  
\]
Let $P_{K}=P_{\mathcal{H}_{K}}, P_{M}=P_{\mathcal{H}_{M}}$ be orthogonal projections onto $\cH_{K}, \cH_{M}$, respectively, and let the \emph{joint concentration} (with respect to $\Lambda_1$ and $\Lambda_2$) $\kappa$ be defined as
\begin{align*}
\kappa &= \sup_{x \in \cH, \enskip y \in \cH_M} \frac{\norm{\mathbbm{1}_{\Lambda_1} \Phi_1(y+ P_K x)}_1 + \norm{\mathbbm{1}_{\Lambda_2} \Phi_2x}_1}{\norm{\Phi_1(y+ P_K x)}_1 + \norm{\Phi_2x}_1}.
\end{align*}
Suppose $(x_{1}^{\star},x_{2}^{\star})$ solves the constrained optimization problem
\begin{align}\label{INPSEP2} (x_{1}^\star ,x_{2}^\star )=\argmin_{(x_{1},x_{2})}\|\Phi_{1}x_{1}\|_{1}+\|\Phi_{2}x_{2}\|_{1} \quad \emph{s.t.} \quad P_{K}(x_{1}+x_{2})=P_{K}(x^{0}).\end{align}  Then $$\|x_{1}^\star -x_{1}^{0}\|_{2}+\|x_{2}^\star -x_{2}^{0}\|_{2}\le \frac{2\delta}{1-2\kappa}.$$
\end{thm}

In Theorem \ref{thm1}, $\cH_{K}$ represents the part of the data which is known, while $\cH_M$ represents the missing part which is to be reconstructed.  In the context of image inpainting, $\mathcal{H}_{K}$ represents the known image pixels, while $\cH_{M}$ represents the unknown pixels to be recovered.  $\Phi_1$ and $\Phi_2$ yield sparse representations of different types of structural information; for example, $\Phi_1$ could be a shearlet frame which sparsely represents cartoon-like images and $\Phi_2$ could sparsely represent texture. If $x^0$ is a superposition of structures which are sparsely represented by the $\Phi_i$, then relatively small index sets $\Lambda_i$ can be chosen yielding both a small $\delta$ and a small $\kappa$.
\begin{proof}
We perform an initial estimate using frame theory and basic $\ell^{p}$ analysis:
\begin{align*}
&\|x_{1}^\star -x_{1}^{0}\|_{2}+\|x_{2}^\star -x_{2}^{0}\|_{2} \\
=& \|\Phi_{1}(x_{1}^\star -x_{1}^{0})\|_{2}+\|\Phi_{2}(x_{2}^\star -x_{2}^{0})\|_{2}\\
\le &\|\Phi_{1}(x_{1}^\star -x_{1}^{0})\|_{1}+\|\Phi_{2}(x_{2}^\star -x_{2}^{0})\|_{1}\\
=&\|\mathbbm{1}_{\Lambda_{1}}\Phi_{1}(x_{1}^\star -x_{1}^{0})\|_{1} + \|\mathbbm{1}_{\Lambda_{2}}\Phi_{2}(x_{2}^\star -x_{2}^{0})\|_{1}+\|\mathbbm{1}_{\Lambda_{1}^{c}}\Phi_{1}(x_{1}^\star -x_{1}^{0})\|_{1} + \|\mathbbm{1}_{\Lambda_{2}^{c}}\Phi_{2}(x_{2}^\star -x_{2}^{0})\|_{1}.
\end{align*}
Noting that $P_{K}(x_{1}^\star -x_{1}^{0})=-P_{K}(x_{2}^\star -x_{2}^{0})$ by the constraint in (\ref{INPSEP2}), $P_M( x_1^0-x_1^\star) + P_K(x_2^\star-x_2^0) = -(x_1^\star -x_1^0)$,
and hence:
\begin{align*}
&\|\mathbbm{1}_{\Lambda_{1}}\Phi_{1}(x_{1}^\star -x_{1}^{0})\|_{1} + \|\mathbbm{1}_{\Lambda_{1}}\Phi_{2}(x_{2}^\star -x_{2}^{0})\|_{1}+\|\mathbbm{1}_{\Lambda_{1}^{c}}\Phi_{1}(x_{1}^\star -x_{1}^{0})\|_{1} + \|\mathbbm{1}_{\Lambda_{2}^{c}}\Phi_{2}(x_{2}^\star -x_{2}^{0})\|_{1}\\
= &\|\mathbbm{1}_{\Lambda_{1}}\Phi_{1}(P_M( x_1^0-x_1^\star) + P_K(x_2^\star-x_2^0))\|_{1} + \|\mathbbm{1}_{\Lambda_{1}}\Phi_{2}(x_{2}^\star -x_{2}^{0})\|_{1}\\
&+\|\mathbbm{1}_{\Lambda_{1}^{c}}\Phi_{1}(P_M( x_1^0-x_1^\star) + P_K(x_2^\star-x_2^0))\|_{1} + \|\mathbbm{1}_{\Lambda_{2}^{c}}\Phi_{2}(x_{2}^\star -x_{2}^{0})\|_{1},
\shortintertext{which, by definition of $\kappa$ and again the constraint in (\ref{INPSEP2}) satisfies}
\le&  \kappa\left(\|\Phi_{1}(P_M( x_1^0-x_1^\star) + P_K(x_2^\star-x_2^0))\|_{1} + \|\Phi_{2}(x_{2}^\star -x_{2}^{0})\|_{1}\right)\\
&+\|\mathbbm{1}_{\Lambda_{1}^{c}}\Phi_{1}(P_M( x_1^0-x_1^\star) + P_K(x_2^\star-x_2^0))\|_{1} + \|\mathbbm{1}_{\Lambda_{2}^{c}}\Phi_{2}(x_{2}^\star -x_{2}^{0})\|_{1}\\
=& \kappa\left(\|\Phi_{1}( x_1^0-x_1^\star)\|_{1} + \|\Phi_{2}(x_{2}^\star -x_{2}^{0})\|_{1}\right)+\|\mathbbm{1}_{\Lambda_{1}^{c}}\Phi_{1}(x_1^0-x_1^\star)\|_{1} + \|\mathbbm{1}_{\Lambda_{2}^{c}}\Phi_{2}(x_{2}^\star -x_{2}^{0})\|_{1}.
\end{align*}

\noindent After some algebraic manipulation and employing $\delta$-sparsity, we see:
\begin{align}
\|\Phi_{1}(x_{1}^\star -x_{1}^{0})\|_{1}+\|\Phi_{2}(x_{2}^\star -x_{2}^{0})\|_{1}&\leq \frac{1}{1-\kappa} \left( \|\mathbbm{1}_{\Lambda_{1}^{c}}\Phi_{1}(x_1^0-x_1^\star)\|_{1} + \|\mathbbm{1}_{\Lambda_{2}^{c}}\Phi_{2}(x_{2}^\star -x_{2}^{0})\|_{1}\right)\nonumber \\
&\leq \frac{1}{1-\kappa} \left( \|\mathbbm{1}_{\Lambda_{1}^{c}}\Phi_{1}x_1^\star\|_{1} + \|\mathbbm{1}_{\Lambda_{2}^{c}}\Phi_{2}x_{2}^\star\|_{1} + \delta \right). \label{eqn:intermezzo2}
\end{align}
So, it suffices to bound these terms involving $\mathbbm{1}_{\Lambda_{i}^{c}}, i=1,2$. We note that $(x_{1}^\star,x_{2}^\star)$ is a minimizer of (\ref{INPSEP2}). Thus, 
\begin{align}
&\|\mathbbm{1}_{\Lambda_{1}^{c}}\Phi_{1}x_1^\star\|_{1} + \|\mathbbm{1}_{\Lambda_{2}^{c}}\Phi_{2}x_{2}^\star\|_{1}\nonumber\\
= & \|\Phi_{1}x_1^\star\|_{1} + \|\Phi_{2}x_{2}^\star\|_{1} -  \|\mathbbm{1}_{\Lambda_{1}}\Phi_{1}x_1^\star\|_{1} - \|\mathbbm{1}_{\Lambda_{2}}\Phi_{2}x_{2}^\star\|_{1}\nonumber\\
\leq&  \|\Phi_{1}x_1^0\|_{1} + \|\Phi_{2}^0 x_{2}^{0}\|_{1} -  \|\mathbbm{1}_{\Lambda_{1}}\Phi_{1}x_1^\star\|_{1} - \|\mathbbm{1}_{\Lambda_{2}}\Phi_{2}x_{2}^\star\|_{1}\nonumber\\
\leq &\|\Phi_{1}x_1^0\|_{1} + \|\Phi_{2}^0 x_{2}^{0}\|_{1} +  \|\mathbbm{1}_{\Lambda_{1}}\Phi_{1}(x_1^\star-x_1^0)\|_{1} -  \|\mathbbm{1}_{\Lambda_{1}}\Phi_{1}x_1^0\|_{1} + \|\mathbbm{1}_{\Lambda_{2}}\Phi_{2}(x_{2}^\star-x_2^0)\|_{1} -  \|\mathbbm{1}_{\Lambda_{2}}\Phi_{2}x_{2}^0\|_{1}\nonumber\\
=& \|\mathbbm{1}_{\Lambda_{1}^c}\Phi_{1}x_1^0\|_{1} +  \|\mathbbm{1}_{\Lambda_{2}^c}\Phi_{2}x_{2}^0\|_{1} +  \|\mathbbm{1}_{\Lambda_{1}}\Phi_{1}(x_1^\star-x_1^0)\|_{1}  + \|\mathbbm{1}_{\Lambda_{2}}\Phi_{2}(x_{2}^\star-x_2^0)\|_{1},\nonumber \\
\shortintertext{which by $\delta$-sparsity is bounded as}
\le& \delta +  \|\mathbbm{1}_{\Lambda_{1}}\Phi_{1}(x_1^\star-x_1^0)\|_{1}  + \|\mathbbm{1}_{\Lambda_{2}}\Phi_{2}(x_{2}^\star-x_2^0)\|_{1}. \label{eqn:part2}
\end{align}
Substituting (\ref{eqn:part2}) into (\ref{eqn:intermezzo2}), we obtain
\begin{align*}
&\|\Phi_{1}(x_{1}^\star -x_{1}^{0})\|_{1}+\|\Phi_{2}(x_{2}^\star -x_{2}^{0})\|_{1}\\
\le & \frac{1}{1-\kappa} \left( \|\mathbbm{1}_{\Lambda_{1}}\Phi_{1}(x_1^\star-x_1^0)\|_{1}  + \|\mathbbm{1}_{\Lambda_{2}}\Phi_{2}(x_{2}^\star-x_2^0)\|_{1}+2\delta \right)\\
\le &\frac{1}{1-\kappa} \left( \kappa\|\Phi_{1}(x_1^\star-x_1^0)\|_{1}  + \|\Phi_{2}(x_{2}^\star-x_2^0)\|_{1}+2\delta \right).
\end{align*}
Thus,
\[
\|\Phi_{1}(x_{1}^\star -x_{1}^{0})\|_{1}+\|\Phi_{2}(x_{2}^\star -x_{2}^{0})\|_{1} \leq \left(1 - \frac{\kappa}{1-\kappa}\right)^{-1} \frac{2 \delta}{1-\kappa} = \frac{2\delta}{1-2\kappa}
.\]\end{proof}
Further theoretical results concerning the joint concentration and applications are the subjects of a forthcoming paper \cite{KMS17}.

\section*{Acknowledgement}
The authors met and began work on this project during the 2016 Hausdorff Trimester Program ``Mathematics of Signal Processing'' and are thus grateful to the organizers and the Hausdorff Institute for Mathematics in Bonn, Germany.

\vspace{\baselineskip}

\bibliography{pamm-KinMur}{}

\end{document}